\documentclass{elsart3-1}

\usepackage{graphicx}
\usepackage{epstopdf}
\usepackage{xcolor}

\usepackage{amssymb}
\usepackage{amsmath}
\usepackage{longtable}

\usepackage[english,french]{babel}
\usepackage[T1]{fontenc}

\newtheorem{theorem}{Theorem}[section]
\newtheorem{conjecture}{Conjecture}[section]
\newtheorem{lemma}[theorem]{Lemma}
\newtheorem{e-proposition}[theorem]{Proposition}
\newtheorem{corollary}[theorem]{Corollary}
\newtheorem{e-definition}[theorem]{Definition\rm}
\newtheorem{remark}{\it Remark\/}
\newtheorem{example}{\it Example\/}

\def\og{\leavevmode\raise.3ex\hbox{$\scriptscriptstyle\langle\!\langle$~}}
\def\fg{\leavevmode\raise.3ex\hbox{~$\!\scriptscriptstyle\,\rangle\!\rangle$}}

\begin{document}
Number theory
\centerline{}

\thanks[label1]{Partially supported by Serbian Ministry of Education and Science, Project 174032.
}
\selectlanguage{english}
\title{The alternative to Mahler measure of a multivariate polynomial}


\selectlanguage{english}
\author{Dragan Stankov},
\ead{dstankov@rgf.bg.ac.rs}

\address{Katedra Matematike RGF-a,
Faculty of Mining and Geology,
University of Belgrade,
Belgrade, \DJ u\v sina 7,
Serbia}
\address{dstankov@rgf.bg.ac.rs}



\selectlanguage{english}
\section*{Abstract}

We introduce the ratio of the number of roots of a polynomial $P_{d}$, less than one in modulus, to its
degree $d$ as an alternative to Mahler measure. We investigate some properties of the alternative. We generalise this definition for
a polynomial in several variables using Cauchy's argument principle.
If a polynomial in two variables do not vanish on the torus we prove the theorem for the alternative which is analogous to the Boyd-Lawton limit formula for Mahler measure. We
determined the exact value of the alternative of $1+x+y$ and $1+x+y+z$. Numerical calculations suggest a conjecture about the exact value
of the alternative of such polynomials having more than three variables.

\selectlanguage{english}
\section{Introduction}

For a non-zero rational function $P \in \mathbb{C}(x_1,\ldots , x_k)$, 
we define the
logarithmic Mahler measure of $P$ to be
$$m(P):=\int_0^1 \int_0^1 \cdots \int_0^1\log|P(e^{2\pi i t_1},e^{2\pi i t_2},\ldots,e^{2\pi i t_k})|dt_1 dt_2\cdots dt_k$$
It is the average value of $\log |P|$ over the unit $k-$torus.
If
\begin{equation}\label{eq:Poly1}
P(x) = a_d \prod_{j=1}^{d}(x - \alpha_j)\in \mathbb{C}_d[x],
\end{equation}
is a polynomial that has only one variable then Jensen’s formula implies
$$m(P)=\int_{0}^1 \log|P(e^{2\pi i t})|dt = \log |a_d| + \sum_{j:|\alpha_j|>1}\log|\alpha_j|. $$
Mahler measure of $P(x)$ is 
$$M(P(x)) :=\exp(m(P))= |a_d|\prod_{j=1}^{d}\max(1, |\alpha_j|).$$
Recall that a cyclotomic (circle dividing)
polynomial is defined as an irreducible factor of
$x^n - 1$, $n \in \mathbb{N}$. Clearly, if $\Phi$
is cyclotomic, then $m(\Phi) = 0$. A well known
Kronecker’s Lemma states that $P \in \mathbb{Z}[x], P \ne 0,
m(P) = 0$ if and only if $P(x) = x^n \prod_{i}\Phi_i(x)$,
where $\Phi_i(x)$ are cyclotomic polynomials.
This also means that polynomials with integer coefficients have logarithmic Mahler
measure greater than or equal to zero.
In 1933 Lehmer asked a question which is still open:
do we have a constant $\delta > 0$ such that for any $P \in \mathbb{Z}[x]$ with
non-zero logarithmic Mahler measure, we must also have $m(P) > \delta$? The smallest known logarithmic Mahler measure $>0$ of a polynomial with integer coefficients is:  

$$m(x^{10} + x^9 - x^7 - x^6 - x^5 - x^4 - x^3 + x + 1) \approx 0.162357612 \ldots .$$

Properties of the logarithmic Mahler measure of a univariate polynomial:

(1) For $P \in Z[x]$, a non zero polynomial, it follows from Kronecker’s Lemma that $M(P) \ge 1$ which implies $m(P) \ge 0$.

(2) For $P, Q \in C[x]$, we have $m(P \cdot Q) = m(P) + m(Q)$, in particular for $c \in C \backslash \{0\}$, we
have $m(cP) = \log |c| + m(P)$.

(3) For a cyclotomic polynomial denoted by $\Phi$, since $M(\Phi) = 1$, we have $m(\Phi) = 0$.

(4) Let P be the product of cyclotomic polynomials, then, we have $m(P) = 0$, and in
addition, for any $n \in N$ we have $m(\pm x^n P) = 0$ .

In general, calculating the Mahler measure of multi-variable polynomials is
much more difficult than the univariate case. 
We have the Boyd-Lawton formula for any rational function
$P \in \mathbb{C}(x_1,\ldots , x_k)$
$$\lim_{ n_2 \to \infty}\cdots \lim_{ n_k \to \infty} m(P(x, x^{n_2},\ldots, x^{n_k})) = m(P(x_1, x_2,\ldots , x_k)),$$
where the $n_j$’s vary independently.

It turns out that for certain polynomials, the Mahler measure is 
in fact a special value of an $L-$function.
In 1981 Smyth proved that
\begin{equation}\label{eq:Smyth1}
m(1 + x + y) =\frac{3\sqrt{3}}{2\pi}L(\chi_{-3}, 2) = L'(\chi_{-3}, -1),
\end{equation}
\begin{equation}\label{eq:Smyth2}
m(1 + x + y + z) = \frac{7}{2\pi^2}\zeta(3) = -14\zeta'(-2).
\end{equation}

Pritsker \cite{Pri} defined a natural areal analog of the Mahler measure and studied its properties. Flammang \cite{F1a1} introduced the absolute S-measure for polynomial \eqref{eq:Poly1} defined by
$$s(P) :=\frac{1}{d}\sum_{j=1}^{d}|{\alpha}_j|,$$ as an analog to Mahler measure and studied its properties. 
For many sequences of polynomials we calculated \cite{Sta1} the limit ratio of the number of roots out of the unit circle to its degree. Each sequence is correlated to a bivariate polynomial having small Mahler measure discovered by Boyd and Mossinghoff \cite{BM}.  

We introduce $c(P)$, our alternative to Mahler measure of a polynomial $P$, as the probability that a randomly chosen zero of $P$ is less than 1 in modulus.
At the beginning of second section of \cite{Pri} 
normalized zero counting measure for polynomial \eqref{eq:Poly1} is defined by
$$\nu_d :=\frac{1}{d}\sum_{j=1}^{d}\delta_{{\alpha}_j},$$ where $\delta_{{\alpha}_j}$ is the unit pointmass at ${\alpha}_j$. 
We can see that $c(P)=\nu_d$ of the open unit disc.

We generalise the definition of $c$ for bivariate polynomials and prove 
\begin{theorem}\label{cha:TrinA}
$$c(a+x+y)=\frac{\arccos(|a|/2)}{\pi}$$
\end{theorem}
The obvious consequence of Theorem \ref{cha:TrinA} is:
\begin{corollary}\label{cha:OneThird}
  $$c(1+x+y)=\frac{1}{3}.$$
\end{corollary}
We generalise the definition of $c$ for multivariate polynomials and prove 
\begin{theorem}\label{cha:OneQuarter}
  $$c(1+x+y+z)=\frac{1}{4}.$$
\end{theorem}
We show that the Boyd-Lawton formula is also valid for our alternative to Mahler measure of a multivariate polynomial $P$ and prove 
its consequence:
\begin{corollary}\label{TwoQuarters}
  $$\lim_{n_2\to\infty} \lim_{n_3\to \infty} c(1+x+x^{n_2}+x^{n_3})=\frac{1}{2},$$
\end{corollary}
where the $n_j$'s vary independently.

\section{The alternative to Mahler measure of a bivariate polynomial}

\begin{theorem}
(Cauchy's argument principle)
If $f(x)$ is a meromorphic function inside and on some closed, simple contour $C$, and $f$ has no zeros or poles on $C$, then
$$\frac{1}{2\pi i}\oint_C \frac{f'(x)}{f(x)}dx=Z-\Pi$$
where $Z$ and $\Pi$ denote respectively the number of zeros and poles of $f(x)$ inside the contour $C$, with each zero and pole counted as many times as its multiplicity and order, respectively, indicate.
\end{theorem}

If $f(x)$ is a polynomial $P(x)$ of degree $d$ and $C$ is the unit circle then we set $x=e^{2\pi i t}$, $dz=2\pi i e^{2\pi i t}$ probability that a randomly chosen zero is inside the unit circle is
$$ c(P):=\frac{1}{d} \int_0^1 \frac{P'(e^{2\pi i t})e^{2\pi i t}}{P(e^{2\pi i t})}dt=\frac{Z}{d}$$ 

Properties of the alternative Mahler measure of univariate polynomials:

(1) For $P \in C[x]$, we have $0 \le c(P) \le 1$.

(2) For $P \in Z[x]$, a non zero polynomial, it follows from Kronecker’s Lemma that $c(P) < 1$.

(3) If $P(x)=a_dx^d+a_{d-1}x^{d-1}+a_0$ and $a_0\ne 0 $ then for the reciprocal polynomial $P^*$ defined $P^*(x)=x^dP(1/x)$ we have 
\begin{equation}\label{eq:recPoly}
c(P^*)=1-c(P).
\end{equation}
(4) For $P, Q \in C[x]$, we have $(\deg(P)+\deg(Q))c(P \cdot Q) = \deg(P)c(P) + \deg(Q)c(Q)$, in particular for $\alpha \in C \backslash \{0\}$, we
have $c(\alpha P) = c(P)$.



Before presenting the following definition, we need to provide an explanation. For a two variable polynomial $P(x,y)$ as in Boyd-Lawton formula we want to determine the limit of the probability $c(P(x,x^n))$ when $n \to \infty$.

Let $P(x,y)$ can be written $$P(x,y)=a_{g}(x)y^{g}+a_{g-1}(x)y^{g-1}+ \cdots + a_0(x).$$ The degree $d$ of $P(x,x^n)$ is
\begin{equation}\label{eq:degree}
d=\deg(a_{g}(x))+ng.
\end{equation}
$$ c(P(x,x^n))=\frac{1}{d} \int_0^1 \frac{\frac{\partial}{\partial x}P(e^{2\pi i t},e^{2\pi i nt})e^{2\pi i t}+n\frac{\partial}{\partial y}P(e^{2\pi i t},e^{2\pi i nt})e^{2\pi i nt}}{P(e^{2\pi i t},e^{2\pi i nt})}dt.$$
We can split the previous integral into two integrals and replace $d$ using \eqref{eq:degree}. The first integral
$$ \frac{1}{\deg(a_{g}(x))+ng} \int_0^1 \frac{\frac{\partial}{\partial x}P(e^{2\pi i t},e^{2\pi i nt})e^{2\pi i t}}{P(e^{2\pi i t},e^{2\pi i nt})}dt \to 0,\;\; n \to \infty.$$
The second integral
$$ \frac{1}{\deg(a_{g}(x))+ng} \int_0^1 \frac{n\frac{\partial}{\partial y}P(e^{2\pi i t},e^{2\pi i nt})e^{2\pi i nt}}{P(e^{2\pi i t},e^{2\pi i nt})}dt \to \frac{1}{g} \int_0^1 \frac{\frac{\partial}{\partial y}P(e^{2\pi i t},e^{2\pi i nt})e^{2\pi i nt}}{P(e^{2\pi i t},e^{2\pi i nt})}dt,$$ as $n$ tends to infinity.
The term $e^{2\pi i t}$ becomes increasingly uncorrelated with $e^{2\pi i nt}$ as $n \to \infty$, so that $nt$ can be replaced with a variable $s$ that vary independently to $t$:
$$ \lim_{n \to \infty} \frac{1}{g} \int_0^1 \frac{\frac{\partial}{\partial y}P(e^{2\pi i t},e^{2\pi i nt})e^{2\pi i nt}}{P(e^{2\pi i t},e^{2\pi i nt})}dt = \frac{1}{g} \int_0^1 \int_0^1 \frac{\frac{\partial}{\partial y}P(e^{2\pi i t},e^{2\pi i s})e^{2\pi i s}}{P(e^{2\pi i t},e^{2\pi i s})}dt ds.$$
It is now clear why we introduce the following
\begin{e-definition}
The alternative Mahler measure for a bivariate polynomial $P(x,y)$, having no zeros on the unit torus, is
$$c(P(x,y)):=\frac{1}{g} \int_0^1 \int_0^1 \frac{\frac{\partial}{\partial y}P(e^{2\pi i t},e^{2\pi i s})e^{2\pi i s}}{P(e^{2\pi i t},e^{2\pi i s})}dt ds.$$
\end{e-definition}


Everest and Ward in \cite{EW} proved their Lemma 3.22 using the Stone-Weierstrass theorem. We present here the lemma as
\begin{lemma}\label{LemmaBoydLawton}
Let $\phi:K^2\rightarrow \mathbb{C}$ be any continuous function. Then
$$\lim_{n\to \infty}\int_{0}^{1}\phi(e^{2\pi i t},e^{2\pi i nt})dt=\int_{0}^{1}\int_{0}^{1}\phi(e^{2\pi i t},e^{2\pi i s})dt\;ds.$$
\end{lemma}

We can prove now that the Boyd-Lawton formula is also valid for our alternative $c$ to Mahler measure.
\begin{corollary} If a polynomial $P(x,y)$ does not vanish on the unit torus then $$\lim_{n \to \infty} c(P(x,x^n))=c(P(x,y)).$$
\end{corollary}

\begin{pf}
We have to prove that
$$ \lim_{n \to \infty} \frac{1}{g} \int_0^1 \frac{\frac{\partial}{\partial y}P(e^{2\pi i t},e^{2\pi i nt})e^{2\pi i nt}}{P(e^{2\pi i t},e^{2\pi i nt})}dt = \frac{1}{g} \int_0^1 \int_0^1 \frac{\frac{\partial}{\partial y}P(e^{2\pi i t},e^{2\pi i s})e^{2\pi i s}}{P(e^{2\pi i t},e^{2\pi i s})}dt ds.$$
Since the integrand is a continuous function it follows that we can use Lemma \ref{LemmaBoydLawton}.
\qed
\end{pf}

In our previous paper \cite{Sta2} we determined the limit ratio 
between number of roots outside the unit circle of
$$P_n(x)=x^n+ax+1$$
and its degree $n$.

\begin{theorem}
The rate $\frac{\nu_{n,a}}{n}$ between the number $\nu_{n,a}$ of roots of the trinomial $x^n+ax+1$, $a\in (0,2]$, which are greater than 1 in modulus, and degree $n$, tends to $\frac{1}{\pi}\arccos(-\frac{a}{2})$, $n\to \infty$.
\end{theorem}
As $P_n(x)=x^n+ax+1=P(x,x^n)$ where $P(x,y)=y+ax+1$ we can use our double integral formula. We can verify that these two values
$$1-\frac{1}{\pi}\arccos(-\frac{a}{2}),$$
$$\int_0^1 \int_0^1 \frac{e^{2\pi i s}}{e^{2\pi i s}+ ae^{2\pi i t}+1}dt ds$$
match to a few decimal places depending of the precision of the definite integral calculator.

\begin{lemma}\label{cha:OneDef}
If $|a|\ne 1$ then
$$\int_{0}^{1}\frac{e^{2\pi i t}}{e^{2\pi i t}+a}dt=
\left\{
	\begin{array}{ll}
		0  & \mbox{if } |a| > 1 \\
		1 & \mbox{if } |a| < 1
	\end{array}
\right.=:1_D(a)$$

\end{lemma}

\begin{pf}
It is the direct consequence of Cauchy's argument principle.
\qed
\end{pf}

\begin{lemma}
If $s\in \left[-\frac{1}{3},\frac{1}{3}\right]$ then $|1+e^{2\pi is}|>1$
\end{lemma}

\begin{pf}
We use the fact that
\begin{equation}\label{eq:e1Cos}
|1+e^{2\pi is}|=2\cos(\pi s).
\end{equation}
\qed
\end{pf}

\begin{pf} of Theorem \ref{cha:TrinA} 
It is convenient to use the obvious fact that $c(a+x+y)=c_x(a+x+y)$. If we: 1. use the definition of $c_x(P(x,y))$, 2. use Lemma \ref{cha:OneDef}, 3. change the variable $2\pi s=\varphi$ in the definite integral, 4. use the determination of $\varphi$ presented on the Figure 1., 5. calculate the definite integral, we obtain the following five steps:
\begin{align*}
  c(a+x+y) & =\int_{0}^{1}\int_{0}^{1}\frac{e^{2\pi i t}}{e^{2\pi i t}+e^{2\pi i s}+a}dtds \\
   & =\int_{0}^{1}1_D(e^{2\pi i s}+a)ds\\
   & =\frac{1}{2\pi}\int_{0}^{2\pi}1_D(e^{\varphi i }+a)d\varphi\\
   & =\frac{1}{2\pi}\int_{\pi-\arccos(|a|/2)}^{\pi +\arccos(|a|/2)}1\cdot d\varphi\\
   & =\frac{\arccos(|a|/2)}{\pi} .\\
\end{align*}
\qed
\end{pf}
We present Smyth's proof \cite{MS,BZ} of \eqref{eq:Smyth1} to illustrate that there is an analogy with our proof of the Theorem \ref{cha:TrinA}.  
\begin{pf} (Smyth 1981) of \eqref{eq:Smyth1}.
We use Jensen's formula: $\int_0^1\log|1+x + e^{2\pi is}|ds=\log^+|1+x|)$
$$m(1 + x + y) = \int_{-1/2}^{1/2}\log^+|1+e^{2\pi it}|dt = \int_{-1/3}^{1/3}\log|1+e^{2\pi it}|dt$$
$$\Re \int_{-1/3}^{1/3}\log(1+e^{2\pi it})dt=\Re \int_{-1/3}^{1/3}\sum_{n=1}^{\infty}\frac{(-1)^{n-1}e^{2\pi int}}{n}dt  $$
$$=\Re \left.\sum_{n=1}^{\infty}\frac{(-1)^{n-1}e^{2\pi int}}{2\pi in^2} \right|_{t=-1/3}^{t=1/3}=\frac{1}{\pi}\sum_{n=1}^{\infty}\frac{(-1)^{n-1}\sin\frac{2\pi n}{3}}{n^2}  $$
It remains to use
$\sin\frac{2\pi n}{3}=\chi_{-3}(n)$ and the multiplicativity of the character $\chi_{-3}$.
\end{pf}
\qed

\begin{figure}[t]
  \includegraphics[width=\textwidth, height=0.38\textheight]{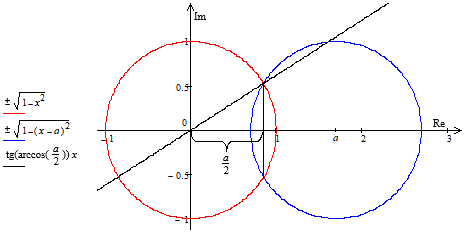}%
\caption{The determination of $\varphi$ such that $e^{\varphi i }+a $ (it is located on the blue circle) is inside of the unit circle (the red circle). It follows that $\varphi \in [\pi -\arccos(|a|/2), \pi +\arccos(|a|/2)]$.}
\label{fig:5}       
\end{figure}


We can show by changing the order of integration in a double integral that switching of the variables in a bivariate polynomial does not effect to the Mahler measure i.e. $m(P(x,y))=m(P(y,x))$.
For the alternative to Mahler measure this is not true. We can demonstrate this with the following 
\begin{example}
If $P(x,y)=xy+x+1$ we can show that $c(P(x,y))\ne c(P(y,x))$: 
\begin{align*}
  c(xy+x+1) & =\lim_{x\to \infty}c(x^{n+1}+x+1) \\
   & =1/3, \\
\end{align*}
using Corollary \ref{cha:OneThird}. On the other hand if we switch the variables and determine $c(P(y,x))$ using Boyd-Lawton formula for $c$ and property \eqref{eq:recPoly} we obtain 
\begin{align*}
  c(yx+y+1) & =\lim_{n\to \infty}c(x^{n+1}+x^n+1) \\
   & =\lim_{n\to \infty}c(x^{n+1}(1+\frac{1}{x}+\frac{1}{x^{n+1}}))\\
   & =\lim_{n\to \infty}c((1+x+x^{n+1})^*)\\
   & =\lim_{n\to \infty}(1-c(x^{n+1}+x+1))\\
   & =1-c(xy+x+1) \\
   & =1-1/3 \\
   & =2/3. \\
\end{align*}
\end{example}

\section{The alternative to Mahler measure of a multivariate polynomial}



Let $P(x_1,x_2,\ldots,x_k)$ can be written $$P(x_1,x_2,\ldots,x_k)=a_{d_j}(x_1,\ldots,x_{j-1},x_{j+1},\ldots,x_k)x_j^{d_j}+$$
$$+a_{d_j-1}(x_1,\ldots,x_{j-1},x_{j+1},\ldots,x_k)x_j^{d_j-1}+ \cdots + a_0(x_1,\ldots,x_{j-1},x_{j+1},\ldots,x_k).$$

\begin{e-definition}
The alternative measure $c_{x_j}$, with respect to $x_j$, or abbreviated $c_j$ for a polynomial in $k$ variables $P(x_1,x_2,\ldots,x_k)$, having no zeros on the unit $k$-torus, is
$$c_j(P(x_1,x_2,\ldots,x_k)):=\frac{1}{d_j} \int_0^1 \int_0^1 \cdots \int_0^1 \frac{\frac{\partial}{\partial x_j}P(e^{2\pi i t_1},e^{2\pi i t_2},\ldots, e^{2\pi i t_k})e^{2\pi i t_j}}{P(e^{2\pi i t_1},e^{2\pi i t_2},\ldots, e^{2\pi i t_k})}dt_1 dt_2\cdots dt_k.$$
\end{e-definition}

\begin{remark}
  If we use $c$ without a subscript it means that the last variable of $P$ should be taken as the subscript so that
  $c(P(x_1,x_2,\ldots,x_k))=c_k(P(x_1,x_2,\ldots,x_k))$ and $c(P(x,y,z))=c_z(P(x,y,z))$.
\end{remark}
In the previous example we showed that  $c(P(x,y))=c_y(P(x,y))=1/3$ and $c(P(y,x))=c_y(P(y,x))=c_x(P(x,y))=2/3$ so that $c_x(P(x,y))\ne c_y(P(x,y)).$


\begin{pf} of Theorem \ref{cha:OneQuarter}.
Again it is convenient to use the obvious fact that $c(1+x+y+z)=c_x(1+x+y+z)$ 1. If we use: the definition of $c_x(1+x+y+z)$, 2. Theorem \ref{cha:TrinA} taking $a=e^{2\pi i u}+1$, 3. Equation \eqref{eq:e1Cos}, 4. the fact that $\arccos \cos(\pi u)=\pi u$, $u\in[0,1/2]$, 5. the calculation of the definite integral, we obtain the following five steps:
\begin{align*} 
  c(1+x+y+z) & =\int_{0}^{1}\int_{0}^{1}\int_{0}^{1}\frac{e^{2\pi i t}}{e^{2\pi i t}+e^{2\pi i s}+e^{2\pi i u}+1}dtdsdu \\
   & =\int_{0}^{1}\frac{\arccos(|e^{2\pi i u}+1|/2)}{\pi}du \\
   & =2\int_{0}^{1/2}\frac{\arccos \cos(\pi u)}{\pi}du\\
   & =2\int_{0}^{1/2}\frac{\pi u}{\pi}du \\
   & =\frac{1}{4} .\\
\end{align*}
\qed
\end{pf}

Numerical calculations of multiple integrals of such functions of $k \le 6$ variables 
$$c_k(1+x_1+x_2+\ldots+x_k)):=\int_0^1 \int_0^1 \cdots \int_0^1 \frac{e^{2\pi i t_k}}{1+e^{2\pi i t_1}+e^{2\pi i t_2}+\ldots+ e^{2\pi i t_k}}dt_1 dt_2\cdots dt_k.$$
suggest us that the following 
conjecture is valid:
\begin{conjecture}\label{cha:OneKth}
$$c(1+x_1+x_2+\cdots+x_k)=\frac{1}{k+1}.$$
\end{conjecture}

\begin{theorem} 
  If $P(x_1,x_2,\ldots,x_k)$ does not vanish on the torus $T^k$ then 
  $$\lim_{n_2\to\infty}\cdots \lim_{n_k\to \infty} c(P(x_1,x_1^{n_2},\ldots,x_1^{n_k}))=c_2(P(x_1,x_2,\ldots,x_k))+\cdots +c_k(P(x_1,x_2,\ldots,x_k)).$$
\end{theorem}

\begin{pf}
Let $P(x_1,x_2,\ldots,x_k)$ can be written $$P(x_1,x_2,\ldots,x_k)=a_{j,d_j}(x_1,\ldots,x_{j-1},x_{j+1},\ldots,x_k)x_j^{d_j}+$$
$$+a_{j,d_j-1}(x_1,\ldots,x_{j-1},x_{j+1},\ldots,x_k)x_j^{d_j-1}+ \cdots + a_{j,0}(x_1,\ldots,x_{j-1},x_{j+1},\ldots,x_k),\;\;j=2,\ldots,k.$$

The degree $d$ of $P(x_1,x_1^{n_2},\ldots,x_1^{n_k})$ is
\begin{equation}\label{eq:degreej}
d=\deg(a_{j,d_j}(x_1,\ldots,x_1^{n_{j-1}},x_1^{n_{j+1}},\ldots,x_1^{n_k}))+n_jd_j.
\end{equation}

$$ c(P(x_1,x_1^{n_2},\ldots,x_1^{n_k}))=$$
$$=\frac{1}{d} \int_0^1 \frac{\frac{\partial}{\partial x_1}P(e^{2\pi i t},e^{2\pi i n_2t},\ldots,e^{2\pi i n_kt})e^{2\pi i t}+ \sum_{j=2}^{k} n_j\frac{\partial}{\partial x_j}P(e^{2\pi i t},e^{2\pi i n_2t},\ldots,e^{2\pi i n_kt})e^{2\pi i n_jt}}  
{P(e^{2\pi i t},e^{2\pi i n_2t},\ldots,e^{2\pi i n_kt})}dt.$$
We can split the previous integral into $k$ integrals and replace $d$ using \eqref{eq:degreej}. The first integral
$$ \frac{1}{\deg(a_{j,d_j}(x_1,\ldots,x_1^{n_{j-1}},x_1^{n_{j+1}},\ldots,x_1^{n_k}))+n_jd_j} \int_0^1 \frac{\frac{\partial}{\partial x_1}P(e^{2\pi i t},e^{2\pi i n_2t},\ldots,e^{2\pi i n_kt})e^{2\pi i t}}{P(e^{2\pi i t},e^{2\pi i n_2t},\ldots,e^{2\pi i n_kt})}dt \to 0,$$
as $n_j \to \infty, j=2,\ldots,k.$ The other integrals

$$\frac{1}{\deg(a_{j,d_j}(x_1,\ldots,x_1^{n_{j-1}},x_1^{n_{j+1}},\ldots,x_1^{n_k}))+n_jd_j} \int_0^1 \frac{ n_j\frac{\partial}{\partial x_j}P(e^{2\pi i t},e^{2\pi i n_2t},\ldots,e^{2\pi i n_kt})e^{2\pi i n_jt}}  
{P(e^{2\pi i t},e^{2\pi i n_2t},\ldots,e^{2\pi i n_kt})}dt.$$ tend to 
$$\frac{1}{d_j} \int_0^1 \int_0^1\cdots \int_0^1 \frac{ \frac{\partial}{\partial x_j}P(e^{2\pi i t_1},e^{2\pi i t_2},\ldots,e^{2\pi i t_k})e^{2\pi i t_j}}{P(e^{2\pi i t_1},e^{2\pi i t_2},\ldots,e^{2\pi i t_k}}dt_1 dt_2\cdots dt_k=c_j,$$ 
$j=2,\ldots,k$ as $n_2,\ldots, n_k \to \infty$ independently, using the following lemma that Boyd proved in Appendix 4 of \cite{Boy1}:
\begin{lemma}
  Suppose $f(x_1,\ldots, x_k)$ is a continuous function on the torus $T^k$, then 
  $$\lim_{n_2\to\infty}\cdots \lim_{n_k\to \infty}\int_{0}^{1}f(e^{2\pi i t},e^{2\pi i n_2t},\ldots,e^{2\pi i n_kt})dt=\int_{0}^{1}\cdots \int_{0}^{1}f(e^{2\pi i t_1},e^{2\pi i t_2},\ldots,e^{2\pi i t_k})dt_1\ldots dt_k.$$
\end{lemma} 
\qed
\end{pf}


\begin{pf} of Corollary \ref{TwoQuarters}
  $$\lim_{n_2\to\infty} \lim_{n_3\to \infty} c(1+x+x^{n_2}+x^{n_3})=c_y(1+x+y+z)+c_z(1+x+y+z)=\frac{1}{4}+\frac{1}{4}=\frac{1}{2}.$$
\qed
\end{pf}
We can show in a similar manner that if Conjecture \ref{cha:OneKth} is true then the following conjecture is also true.
\begin{conjecture}
  $$\lim_{n_2\to\infty} \lim_{n_3\to \infty}\cdots \lim_{n_k\to \infty} c(1+x+x^{n_2}+x^{n_3}+\cdots+x^{n_k})=\frac{k-1}{k+1},$$
\end{conjecture}
where the $n_j$'s vary independently.


\begin{thebibliography}{00}
\bibitem{BCFJ} P. Borwein, S. Choi, R. Ferguson, and J. Jankauskas, On Littlewood polynomials with prescribed number of zeros inside the unit disk, Canad. J. of Math. 67 (2015) 507--526.

\bibitem{BEFL} P. Borwein, T. Erd\'{e}lyi, R. Ferguson, and R. Lockhart, On the zeros of cosine polynomials:
solution to a problem of Littlewood, Ann. Math. Ann. (2) 167 (3) (2008) 1109–1117.


\bibitem{Boy1} D. W. Boyd, Speculations concerning the range of Mahler's measure. Canad. Math.
Bull. 24 (4) (1981) 453 -- 469.

\bibitem{BM} D. W. Boyd, M. J. Mossinghoff, Small limit points of Mahler’s measure. Experiment. Math. 14 (2005), No. 4, 403--414

\bibitem{BZ} F. Brunault, W. Zudilin. Many Variations of Mahler Measures. A Lasting Symphony. Cambridge University Press, (2020). 

\bibitem{BGMP} Brunault, François; Guilloux, Antonin; Mehrabdollahei, Mahya; Pengo, Riccardo. Limits of Mahler measures in multiple variables. Annales de l'Institut Fourier, Volume 74 (2024) no. 4, pp. 1407--1450. 


\bibitem{Dru} P. Drungilas, Unimodular roots of reciprocal Littlewood polynomials, J. Korean Math. Soc. 45 (3)(2008) 835--840.


\bibitem{Dub2} A. Dubickas, (2023). Every Salem number is a difference of two Pisot numbers. Proc. Edinb. Math. Soc., 66(3), 862--867. 



\bibitem{EW} G. Everest, T. Ward, Heights of Polynomials and Entropy in Algebraic Dynamics, Springer-Verlag London Ltd., London, (1999).

\bibitem{F1a1} V. Flammang. The S-measure for algebraic integers having all their conjugates in a sector. Rocky
Mountain Journal of Mathematics, Rocky Mountain Mathematics Consortium, 50 (4), (2020) 1313--1321.


\bibitem{FV} V. Flammang, P. Voutier. Properties of trinomials of height at least 2, \emph{Rocky Mountain J. Math.}, 2022, 52(2), 507 -- 518.

\bibitem{GL} J. Gu, M. Lalín, The Mahler measure of a three-variable family and an application to the Boyd-Lawton formula, Res. Number Theory 7 (2021), no. 1, article
no. 13 (23 pages).

\bibitem{GVG} C. Guichard and J.-L. Verger-Gaugry, On Salem numbers, expansive polynomials and Stieltjes continued fractions, J. Th\'{e}orie Nombres Bordeaux 27,
(2015), 769--804.

\bibitem{LN} M. Lal{\'i}n, S.S. Nair, An invariant property of Mahler measure, Bulletin of the London Mathematical Society, vol. 55 (3), (2023), 1129--1142

\bibitem{MS} J. McKee, C. Smyth, Around the Unit Circle, Springer International Publishing, London, (2021). ISBN: 9783030800307

\bibitem{Muk} K. Mukunda, Littlewood Pisot numbers, J. Number Theory 117 (1) (2006) 106--121.

\bibitem{Pri} Pritsker, Igor. (2008). An areal analog of Mahler's measure. Illinois Journal of Mathematics. 52. 347--363. 




\bibitem{Sta1} D. Stankov. The number of nonunimodular roots of a reciprocal polynomial. Comptes Rendus. Math\'{e}matique, Volume 361 (2023), p. 423--435, https://doi.org/10.5802/crmath.422

\bibitem{Sta2} D. Stankov. The Boyd's conjecture. https://arxiv.org/abs/arXiv:1401.1688v2, March
2014.








\end{thebibliography}
\end{document}